\documentclass{pja00}
\usepackage{amssymb}
\usepackage{amscd}
\usepackage{latexsym}
\usepackage{eucal}
\usepackage{enumerate}
\usepackage[all]{xy}
\usepackage{graphics}

\newtheorem{theorem}{Theorem}

\newtheorem{definition}[theorem]{Definition}
\newtheorem{fact}[theorem]{Fact}
\newtheorem{remark}[theorem]{Remark}
\newtheorem{question}[theorem]{Questions}
\newtheorem{problems}[theorem]{Problems}
\newtheorem{example}[theorem]{Example}
\newtheorem{setting}[theorem]{Setting}

\numberwithin{equation}{section}%

\runninghead{Spectral analysis on pseudo-Riemannian locally symmetric spaces}
\title{Spectral analysis on pseudo-Riemannian locally symmetric spaces}

\Author{1}{Fanny}{Kassel}

\Author{2}{Toshiyuki}{Kobayashi}

\affiliation{1}{CNRS and IHES,
Laboratoire Alexander Grothendieck,
35 route de Chartres, 
91440 Bures-sur-Yvette, France.
Supported by the European Research Council under the European Union's Horizon 2020 research and innovation programme (ERC starting grant DiGGeS, grant agreement No 715982).}

\affiliation{2}{Graduate School of Mathematical Sciences,
The University of Tokyo,
3-8-1 Komaba, Tokyo, 153-8914, Japan,
and 
Kavli Institute for the Physics and Mathematics of the Universe (WPI).
Partially supported  by Kakenhi JP18H03669.
}

\begin{abstract}
We summarize recent results initiating spectral analysis on pseudo-Riemannian locally symmetric spaces $\Gamma \backslash G/H$, beyond the classical setting where $H$ is compact (e.g.\ theory of automorphic forms for arithmetic $\Gamma$) or $\Gamma$ is trivial (e.g.\ Plancherel-type formula for semisimple symmetric spaces).  
\end{abstract}


\Subject[2020]{22E40; 22E46, 58J50, 11F72, 53C35}

\date{\today}

\begin{document}

\maketitle%

\section{Introduction} \label{sec:Intro}
\setcounter{subsection}{0}

A {\it{pseudo-Riemannian manifold}} 
 is a smooth manifold $M$
 equipped with a smooth, nondegenerate symmetric
 bilinear tensor $g$ of signature $(p,q)$.  
It is called Riemannian 
 if $q=0$, 
 and Lorenzian if $q=1$.  
As in the Riemannian case, 
 the metric $g$ induces a Radon measure on $M$
 and a second-order
 differential operator
\[
  \square_M=\operatorname{div} \operatorname{grad}
\]
 called the \emph{Laplacian}.
It is a symmetric operator
on the Hilbert space $L^2(X_\Gamma)$.  
The Laplacian $\square_M$ is not an elliptic differential operator
 if $p,q>0$.  

A {\it{semisimple symmetric space}} $X$ is a homogeneous space
 $G/H$
 where $G$ is a semisimple Lie group 
 and $H$ an open subgroup 
 of the group of fixed points of $G$
 under some involutive automorphism.  
The manifold $X$ carries a $G$-invariant pseudo-Riemannian metric induced by the Killing form
 of the Lie algebra ${\mathfrak {g}}$ of $G$.  
The group $G$ acts on $X$ by isometries, 
 and the ${\mathbb{C}}$-algebra
 ${\mathbb{D}}_G(X)$ 
 of $G$-invariant  differential operators
 on $X$
 is commutative.  

In this note we consider quotients $X_{\Gamma}=\Gamma \backslash X$ of a semisimple symmetric space $X=G/H$ by discrete subgroups $\Gamma$ of $G$ acting properly discontinuously and freely on~$X$ (\lq\lq{discontinuous groups for $X$}\rq\rq).\,Such quotients are called \emph{pseudo-Riemannian locally symmetric spaces}.
They are complete $(G,X)$-manifolds in the sense of Ehresmann and Thurston, and they inherit a pseudo-Riemannian structure from~$X$.  
Any $G$-invariant differential operator $D$ on~$X$ induces a differential operator
 $D_{\Gamma}$
 on $X_{\Gamma}$ via the covering map 
 $p_{\Gamma} \colon X \to X_{\Gamma}$.  
For instance, 
 the Laplacian $\square_X$ on $X$ is $G$-invariant, 
 and $(\square_X)_{\Gamma}=\square_{X_{\Gamma}}$.  
We think of 
\[
     \mathcal{P}:=\{D_{\Gamma}: D \in {\mathbb{D}}_G(X)\}
\]
 as the set of
 \lq\lq{intrinsic differential operators}\rq\rq\ 
 on the locally symmetric space $X_{\Gamma}$.  
It is a subalgebra
 of the ${\mathbb{C}}$-algebra
 ${\mathbb{D}}(X_{\Gamma})$ 
 of differential operators on $X_{\Gamma}$:
\begin{equation}
\label{eqn:DGamma}
  {\mathbb{D}}_G(X) 
  \overset \sim \to
  \mathcal{P} \subset {\mathbb{D}}(X_{\Gamma}), 
  \quad
  D \mapsto D_{\Gamma}.  
\end{equation}

For a ${\mathbb{C}}$-algebra homomorphism
 $\lambda \colon {\mathbb{D}}_G(X) \to\nolinebreak {\mathbb{C}}$, 
 we denote by $C^{\infty}(X_{\Gamma};{\mathcal{M}}_{\lambda})$
 the space of smooth functions $f$ on~$X_{\Gamma}$
({\it{joint eigenfunctions}})
 satisfying the following system of partial differential equations:
\begin{equation*}
({\mathcal{M}}_{\lambda})
\qquad\quad
  D_{\Gamma} f = \lambda(D) f
\quad
\text{for all $D \in {\mathbb{D}}_G(X)$.  }
\end{equation*}
Let  $L^2(X_{\Gamma};{\mathcal{M}}_{\lambda})$
be the space of square-integrable functions on~$X_{\Gamma}$
 satisfying $({\mathcal{M}}_{\lambda})$
 in the weak sense.
It is a closed subspace of the Hilbert space 
 $L^2(X_{\Gamma})$.
 We are interested in the following problems.

\begin{problems} 
For intrinsic differential operators
on $X_\Gamma=\Gamma\backslash G/H$,
\label{prob:1}
\begin{enumerate}
\item[{\rm{(1)}}]
construct joint eigenfunctions on $X_{\Gamma}$;
\item[{\rm{(2)}}]
find a spectral theory on $L^2(X_{\Gamma})$.  
\end{enumerate}
\end{problems}
In the classical setting where $H$ is a maximal compact subgroup $K$
 of $G$, i.e.\ $X_{\Gamma}$ is a {\it{Riemannian}} locally symmetric space, a rich and deep theory has been developed
over several decades, in particular,
in connection with automorphic forms
when $\Gamma$ is arithmetic.
For compact $H$, the spectral decomposition 
 of $L^2(X_{\Gamma})$
 is closely related to a disintegration 
 of the regular representation of $G$ on $L^2({\Gamma}\backslash G)$:
\begin{equation}
\label{eqn:direct}
 L^2({\Gamma}\backslash G)
 \simeq
  \int_{\widehat G}^{\oplus} m_{\Gamma}(\pi) \, \pi \, \mathrm{d} \sigma(\pi), 
\end{equation}
where $\mathrm{d} \sigma$ is a Borel measure on the unitary dual $\widehat G$
 and $m_{\Gamma} \colon \widehat G \to {\mathbb{N}} \cup \{\infty\}$
 a measurable function called {\it{multiplicity}}.  
There is a natural isomorphism
\begin{equation}
\label{eqn:Kinv}
  L^2(X_{\Gamma}) \simeq L^2(\Gamma \backslash G)^H
\end{equation}
and the Hilbert space $L^2(X_{\Gamma})$
 is decomposed as 
\begin{equation}
\label{eqn:absP}
  L^2(X_{\Gamma}) \simeq \int_{(\widehat G)_H} m_{\Gamma}(\pi) \, \pi^H \, \mathrm{d} \sigma(\pi), 
\end{equation}
 where $\pi^H$ denotes the space
 of $H$-invariant vectors 
 in the representation space of $\pi$
 and 
\[
   (\widehat G)_H: = \big\{ \pi \in \widehat G: \pi^H \neq \{0\} \big\}.  
\]
Since the center ${\mathfrak{Z}}({\mathfrak{g}}_{\mathbb{C}})$
 of the enveloping algebra 
 $U({\mathfrak{g}}_{\mathbb{C}})$ acts on the space of smooth vectors of $\pi$
 as scalars for every $\pi \in \widehat G$, 
 the decomposition \eqref{eqn:absP} respects the actions
 of ${\mathbb{D}}_G(X)$ and ${\mathfrak{Z}}({\mathfrak{g}}_{\mathbb{C}})$
 via the natural ${\mathbb{C}}$-algebra homomorphism
$
  {\mathrm{d}}\ell \colon{\mathfrak{Z}}({\mathfrak{g}}_{\mathbb{C}}) \to {\mathbb{D}}_G(X). 
$
This homomorphism is surjective e.g.\ if $G$ is a classical group.  

The situation changes drastically 
 beyond the aforementioned classical setting, 
 namely, 
 when $H$ is not compact anymore.  
New difficulties include:
\begin{enumerate}
\item[(1)]
(Representation theory)\enspace
By the ergodicity theorem of Howe--Moore \cite{howemoore}, if $H$ is noncompact, then $L^2(\Gamma \backslash G)^H =\{0\}$, and so \eqref{eqn:Kinv} fails:
\begin{equation}
\label{eqn:Hinv}
  L^2(X_{\Gamma}) \not\simeq L^2(\Gamma \backslash G)^H
\end{equation}
 and the irreducible decomposition \eqref{eqn:direct} of the regular representation $L^2(\Gamma \backslash G)$
 of~$G$ does not yield a spectral decomposition of $L^2(X_{\Gamma})$.  
\item[(2)]
(Analysis)\enspace
In contrast to the usual Riemannian case (see \cite{str83}),
the Laplacian $\square_{X_{\Gamma}}$ is not elliptic anymore,
and thus even the following subproblems of Problem \ref{prob:1}.(2)
 are open in general 
 for $X_{\Gamma}= \Gamma \backslash G/H$
 with $H$ noncompact.  
\end{enumerate}

\begin{question}
\label{problems}
\begin{enumerate}
\item[{\rm{(1)}}] Does the Laplacian $\square_{X_{\Gamma}}$, defined on $C^{\infty}_c(X_{\Gamma})$, extend to a self-adjoint operator on $L^2(X_{\Gamma})$?
\item[{\rm{(2)}}] Does $L^2(X_{\Gamma};\mathcal{M}_{\lambda})$ contain real analytic functions as a dense subspace?
\item[{\rm{(3)}}] 
Does $L^2(X_{\Gamma})$ decompose discretely into a sum
 of subspaces $L^2(X_{\Gamma};\mathcal{M}_{\lambda})$ 
 when $X_{\Gamma}$ is compact?
\end{enumerate}
\end{question}

\section{Standard quotients} \label{subsec:CK}

We observe that a discrete group of isometries on a pseudo-Riemannian manifold $X$
 does not always act properly discontinuously 
 on $X$, 
 and the quotient space 
 $X_{\Gamma} =\Gamma \backslash X$ is not necessarily 
 Hausdorff.  
In fact, 
 some semisimple symmetric spaces $X$ do not admit infinite discontinuous groups
 of isometries 
 (Calabi--Markus phenomenon 
 \cite{CM, kob89}), 
and thus it is not obvious a priori
 whether there are interesting examples
 of pseudo-Riemannian locally symmetric spaces $X_{\Gamma}$ beyond the classical Riemannian case.

Fortunately,
 there exist semisimple symmetric spaces $X=G/H$
 admitting \lq\lq{large}\rq\rq\
 discontinuous groups $\Gamma$
 such that $X_{\Gamma}$ is compact
 or of finite volume.  
Let us recall a useful idea for finding such $X$ and~$\Gamma$.
Suppose a Lie subgroup $L$ of $G$ 
 acts properly on $X$.  
Then the action of any discrete subgroup $\Gamma$ of $L$ on $X$
 is automatically properly discontinuous, 
 and this action is free
 whenever $\Gamma$ is torsion-free.  
Moreover, if $L$ acts cocompactly (e.g.\ transitively) on~$X$,
 then $\operatorname{vol}(X_{\Gamma})<+\infty$
 if and only if $\operatorname{vol}(\Gamma\backslash L)<+\infty$.  

\begin{definition}
[Standard quotient $X_{\Gamma}$]
\label{def:standard}
{\rm{A quotient $X_{\Gamma}=\Gamma \backslash X$
 of $X=G/H$ by a discrete subgroup of $G$
 is called {\it{standard}}
 if $\Gamma$ is contained in a reductive subgroup $L$
 of $G$
 acting properly on $X$.  
}}
\end{definition}
A criterion on triples $(G,L,H)$ of reductive Lie groups for $L$ to act properly on $X=G/H$
 was established in \cite{kob89}, 
 and a list of irreducible symmetric spaces $G/H$ admitting proper
 and cocompact actions
 of reductive subgroups $L$
 was given in \cite{ky05}.  
Recently, 
 Tojo \cite{tojo} announced 
 that the list in \cite{ky05}
 exhausts all 
 such triples $(L, G, H)$ 
 with $L$ maximal.

\section{Construction of discrete spectrum} \label{sec:disc}

Let $X=G/H$ be a semisimple symmetric space.
Let ${\mathfrak{j}}$ be a maximal semisimple abelian subspace
 in the orthogonal complement of ${\mathfrak{h}}$ in ${\mathfrak{g}}$
 with respect to the Killing form, 
 and $W$ the Weyl group
 for the root system $\Sigma({\mathfrak{g}}_{\mathbb{C}}, {\mathfrak{j}}_{\mathbb{C}})$.  
The Harish-Chandra isomorphism
$\Psi\colon S({\mathfrak j}_{\mathbb C})^W
\overset \sim \to \mathbb D_G(X)$ (see \cite{hel00}) induces a
 bijection 
\begin{equation}
\label{eqn:HCisom}
\Psi^{\ast} \colon \operatorname{Hom}_{\mathbb{C}\text{-alg}}({\mathbb{D}}_G(X), {\mathbb{C}}) \overset{\sim}{\longrightarrow} {\mathfrak{j}}_{\mathbb{C}}^{\ast}/W.  
\end{equation}
The dimension of ${\mathfrak{j}}$ is called 
the \emph{rank} of the symmetric space $X=G/H$.  
Let $K$ be a maximal compact subgroup of $G$
such that $H \cap K$ is a maximal compact subgroup of $H$.  
Assume 
 that $G$ is connected 
 without compact factor and that the following rank condition is satisfied:
\begin{equation}
\label{eqn:rank}
\operatorname{rank} G/H
=
\operatorname{rank} K/(H\cap K).  
\end{equation}
Then we can take ${\mathfrak{j}}$ as a subspace of ${\mathfrak{k}}$.  
We fix compatible positive systems 
 $\Sigma^+({\mathfrak{g}}_{\mathbb{C}}, {\mathfrak{j}}_{\mathbb{C}})$
 and $\Sigma^+({\mathfrak{k}}_{\mathbb{C}}, {\mathfrak{j}}_{\mathbb{C}})$, 
 denote by $\rho$ and $\rho_c$ the corresponding half sums of positive roots
 counted with multiplicities, and set 
\begin{equation*}
\Lambda:=
2 \rho_c - \rho+
{\mathbb{Z}}\operatorname{-span}
\big\{\text{highest weights of $(\widehat {K})_{H \cap K}$}\big\}.  
\end{equation*}

For $C \ge 0$, 
 we consider the countable set
\[
  \Lambda_C:=
  \big\{\lambda \in {\Lambda}:
  \langle \lambda, \alpha \rangle >C
\ \text{for all }\alpha \in \Sigma^+({\mathfrak{g}}_{\mathbb{C}}, {\mathfrak{j}}_{\mathbb{C}})\big\}.  
\]

\begin{fact}
[Flensted-Jensen {\cite{fle80}}]
\label{fact:FJ}
If the rank condition \eqref{eqn:rank} holds,
 then there exists $C>0$ such that
\[
  L^2(X;{\mathcal{M}}_{\lambda})\ne \{0\}
\quad
\text{for all }\lambda \in \Lambda_C.  
\]
\end{fact}

In fact one can take $C=0$ \cite{mo84}.
We now turn to locally symmetric spaces $X_\Gamma$:

\begin{theorem}
[{\cite{kk11}}, {\cite[Th.\,1.5]{kk16}}]
\label{thm:Pseries}
Under the rank condition \eqref{eqn:rank},
for any standard quotient $X_{\Gamma}$
 with $\Gamma$ torsion-free, 
 there exists $C_{\Gamma}>0$ such that 
\[
  L^2(X;{\mathcal{M}}_{\lambda}) \ne \{0\}
\quad
 \text{for all } \lambda \in \Lambda_{C_{\Gamma}}.  
\]
\end{theorem}
Thus the \emph{discrete spectrum} $\mathrm{Spec}_d(X_{\Gamma})$, which is by definition the set of $\lambda \in \operatorname{Hom}_{\mathbb{C}\text{-alg}}({\mathbb{D}}_G(X), {\mathbb{C}})$ such that $L^2(X;{\mathcal{M}}_{\lambda}) \ne \{0\}$, is infinite.
 
Theorem~\ref{thm:Pseries} applied to $(G \times \{1\}, G \times G, \operatorname{Diag}G)$ instead of $(L,G,H)$ (group manifold case) implies:

\begin{example}
{\rm{
Suppose $\operatorname{rank}G=\operatorname{rank} K$.  
For {\it{any}} torsion-free discrete subgroup $\Gamma$ 
 and {\it{any}} discrete series representation $\pi_{\lambda}$
 of~$G$ with sufficiently regular Harish-Chandra parameter $\lambda$,
\begin{equation}
\label{eqn:autoform}
   \operatorname{Hom}_G(\pi_{\lambda}, L^2(\Gamma \backslash G))
   \ne \{0\}.  
\end{equation}
This sharpens and generalizes the known results
 asserting that 
 if $\Gamma$ is an {\it{arithmetic}} subgroup of $G$,
 then \eqref{eqn:autoform} holds
 after replacing $\Gamma$ 
 by a finite-index subgroup $\Gamma'$
 (possibly depending on $\pi_{\lambda}$), 
 see Borel--Wallach \cite{bw00}, 
 Clozel \cite{clo86}, 
DeGeorge--Wallach \cite{dw78}, 
Kazhdan \cite{kaz77}, 
 Rohlfs--Speh \cite{rs87},
and 
 Savin \cite{sav89}.
 }}
\end{example}

\begin{remark} \label{rem:Pseries-rank}
{\rm{
{\rm{(1)}}\enspace
Theorem~\ref{thm:Pseries} extends to a more general setting
 where $X_{\Gamma}$ is not necessarily standard: namely,
 the conclusion still holds
 as soon as the action of $\Gamma$ on $X$ satisfies
 a strong properness condition called {\it{sharpness}} \cite[Th.\,3.8]{kk16}.  
\par\noindent
{\rm{(2)}}\enspace
The rank condition \eqref{eqn:rank} is necessary
for $\mathrm{Spec}_d(X)$ to be nonempty
 (see Matsuki--Oshima \cite{mo84}),
 in which case Fact \ref{fact:FJ} applies.
 On the other hand,
 $\mathrm{Spec}_d(X_{\Gamma})$ may be nonempty
 even if \eqref{eqn:rank} fails.  
 This leads us to the notion of discrete spectrum 
 of type {\textbf {I}} and {\textbf {II}},
 see Definition \ref{def:type} below.  
}}
\end{remark}

\section{Spectral decomposition of $L^2(X_{\Gamma})$} \label{sec:decomp}

In this section,
 we discuss spectral decomposition 
 on standard quotients $X_{\Gamma}$.  
We do not impose the rank condition \eqref{eqn:rank}, 
 but require
 that $L_{\mathbb{C}}$ act spherically on $X_{\mathbb{C}}$, 
 i.e.\ a Borel subgroup of $L_{\mathbb{C}}$ has an open orbit in $X_{\mathbb{C}}$.  
To be precise,
 our setting is as follows:
\begin{setting}
\label{setting:main}
We consider a symmetric space $X=G/H$ with $G$ noncompact and simple, a reductive subgroup $L$ of~$G$ acting properly on~$X$ such that $X_{\mathbb{C}}=G_{\mathbb{C}}/H_{\mathbb{C}}$
 is $L_{\mathbb{C}}$-spherical, and a torsion-free discrete subgroup $\Gamma$ of~$L$.  
\end{setting}

For compact $H$, we can take $L=G$. 
However, our main interest is for {\it{noncompact}}~$H$, in which case $L\neq G$ in the setting~\ref{setting:main}.  

In Theorems~\ref{thm:Fourier} and \ref{thm:selfadj} below, 
 we allow the case where $\operatorname{vol}(X_{\Gamma})=+\infty$.  

\begin{theorem}[Spectral decomposition] \label{thm:Fourier}
In the setting~\ref{setting:main},
 there exist a measure $\mathrm{d}\mu$ on $\mathrm{Hom}:=\mathrm{Hom}_{\mathbb{C}\text{-}\mathrm{alg}}(\mathbb{D}_G(X),\mathbb{C})$
 and a measurable family $(\mathcal{F}_{\lambda})_{\lambda\in\mathrm{Hom}}$ of linear maps, with
 $$
  \mathcal{F}_{\lambda} : C_c^{\infty}(X_{\Gamma}) \longrightarrow C^{\infty}(X_{\Gamma};\mathcal{M}_{\lambda}),
$$
such that any $f\in C_c^{\infty}(X_{\Gamma})$ can be expanded into joint eigenfunctions on~$X_{\Gamma}$ as
\begin{equation}
\label{eqn:Fourier}
f = \int_{\mathrm{Hom}} \mathcal{F}_{\lambda} f \ \mathrm{d}\mu(\lambda),
\end{equation}
with a Parseval--Plancherel type formula
\[
  \| f \|_{L^2(X_{\Gamma})}^2
  =
  \int_{\mathrm{Hom}} \| {\mathcal{F}}_{\lambda} f\|_{L^2(X_{\Gamma})}^2
\ \mathrm{d} \mu(\lambda).
\]
\end{theorem}

The measure $\mathrm{d} \mu$ can be described 
 via a \lq\lq{transfer map}\rq\rq\
 discussed in Section \ref{sec:Transfer}, 
 see \eqref{eqn:dmu}.  
In particular, we see that
 \eqref{eqn:Fourier} is a discrete sum if $X_{\Gamma}$ is compact, 
 answering Question \ref{problems}.(3) in our setting.
The proof of Theorem~\ref{thm:Fourier} gives an answer
 to Questions \ref{problems}.(1)--(2):

\begin{theorem}
\label{thm:selfadj}
In the setting~\ref{setting:main}, 
\par\noindent
{\rm{(1)}}\enspace
the pseudo-Riemannian Laplacian $\square_{X_{\Gamma}}$ defined on $C_c^{\infty}(X_{\Gamma})$ is 
essentially self-adjoint on $L^2(X_{\Gamma})$;
\par\noindent
{\rm{(2)}}\enspace
any $L^2$-eigenfunction of the Laplacian $\square_{X_{\Gamma}}$
 can be approximated by
 real analytic $L^2$-eigenfunctions.  
\end{theorem}

\begin{theorem} \label{thm:mainII}
In the setting~\ref{setting:main}, 
the discrete spectrum $\mathrm{Spec}_d(X_{\Gamma})$ is infinite whenever $\Gamma$ is cocompact or arithmetic in the subgroup~$L$.  
\end{theorem}

Let ${\mathcal{D}}'(X)$ be the space of distributions
on $X$, endowed with its standard topology.
Let
$
  p_{\Gamma}^{\ast} \colon L^2(X_{\Gamma}) \to\nolinebreak {\mathcal{D}}'(X)
$
 be the pull-back by the projection 
 $p_{\Gamma}\colon X \to\nolinebreak X_{\Gamma}$.  
For $\lambda\in\mathrm{Spec}_d(X_{\Gamma})$, we denote by
 $L^2(X_{\Gamma};{\mathcal{M}}_{\lambda})_{\bf{I}}$
 the preimage under $p_{\Gamma}^{\ast}$
 of the closure in ${\mathcal{D}}'(X)$
 of $L^2(X_{\Gamma};{\mathcal{M}}_{\lambda})$,
 and by $L^2(X_{\Gamma};{\mathcal{M}}_{\lambda})_{\bf{II}}$
 its orthogonal complement 
 in $L^2(X_{\Gamma};{\mathcal{M}}_{\lambda})$.  

\begin{definition}
\label{def:type}
{\rm{
For $i={\bf{I}}$ or~${\bf{II}}$, the \emph{discrete spectrum of type~$i$} of~$X_{\Gamma}$ is the subset $\mathrm{Spec}_d(X_{\Gamma})_i$ of $\mathrm{Spec}_d(X_{\Gamma})$ consisting of those elements $\lambda$ such that $L^2(X_{\Gamma};{\mathcal{M}}_{\lambda})_i \ne \{0\}$.
}}
\end{definition}

By construction, $\mathrm{Spec}_d(X_{\Gamma})_{\bf{I}}$ is contained in $\mathrm{Spec}_d(X)$, hence it is nonempty only if \eqref{eqn:rank} holds (Remark~\ref{rem:Pseries-rank}.(2)); in this case $\mathrm{Spec}_d(X_{\Gamma})_{\bf{I}}$ is actually infinite for standard $X_{\Gamma}$ by Theorem~\ref{thm:Pseries}.
On the other hand, Theorem~\ref{thm:mainII} has 
the following refinement.

\begin{theorem}
In the setting \ref{setting:main}, 
 $\operatorname{Spec}_d(X_{\Gamma})_{\bf{II}}$ is infinite
 whenever $\Gamma$ is cocompact or arithmetic in $L$.  
\end{theorem}

\begin{example}
{\rm{
Let $M$ be a 3-dimensional compact standard anti-de Sitter manifold.  
Then both $\operatorname{Spec}_d(X_{\Gamma})_{\bf{I}}$ and 
 $\operatorname{Spec}_d(X_{\Gamma})_{\bf{II}}$ are infinite, 
 and 
\[
   \operatorname{Spec}_d(X_{\Gamma})_{\bf{I}} \subset [0,+\infty), 
\quad
  \operatorname{Spec}_d(X_{\Gamma})_{\bf{II}} \subset (-\infty,0].  
\]
}}
\end{example}

\section{Transfer maps} \label{sec:Transfer}

In Section \ref{sec:Intro}
 we considered spectral analysis on locally symmetric spaces $X_{\Gamma}$ through the algebra $\mathcal{P}$ $(\simeq {\mathbb{D}}_G(X))$
of intrinsic differential operators on $X_{\Gamma}$.  
For standard quotients $X_{\Gamma}$ 
 with $\Gamma \subset L$, another ${\mathbb{C}}$-algebra $\mathcal{Q}$ of differential operators on~$X_{\Gamma}$ is obtained from the center ${\mathfrak{Z}}({\mathfrak{l}}_{\mathbb{C}})$
 of the enveloping algebra
 $U({\mathfrak{l}}_{\mathbb{C}})$: indeed, ${\mathfrak{Z}}({\mathfrak{l}}_{\mathbb{C}})$ acts on smooth functions on~$X$ by differentiation, yielding a $\mathbb{C}$-algebra of
 $L$-invariant differential operators on~$X$,
 hence a $\mathbb{C}$-algebra of differential operators on $X_{\Gamma}=\Gamma\backslash X$ since $\Gamma\subset L$.
In general, 
 there is no inclusion relation 
 between $\mathcal{P}$ and $\mathcal{Q}$. 
In order to compare the roles of $\mathcal{P}$ and $\mathcal{Q}$,
we highlight a natural homomorphism ${\mathfrak{Z}}({\mathfrak{g}}_{\mathbb{C}}) \to \mathcal{P}$
and a surjective one $\mathrm{d}\ell \colon {\mathfrak{Z}}({\mathfrak{l}}_{\mathbb{C}}) \to\nolinebreak \mathcal{Q}$.
Loosely speaking, 
 the algebras ${\mathfrak{Z}}({\mathfrak{g}}_{\mathbb{C}})$
 and ${\mathfrak{Z}}({\mathfrak{l}}_{\mathbb{C}})$ separate
 irreducible representations of the groups $G$ and $L$, 
 respectively, 
 hence it is important to understand
 how irreducible representations of $G$
 behave when restricted to the subgroup $L$
 ({\it{branching problem}}) in order to utilize 
 the algebra $\mathcal{Q}$
 for the spectral analysis on $X_{\Gamma}$
 via the algebra $\mathcal{P}$
 (see \cite{kob09, kob17}).  
We shall return to this point in Theorem~\ref{thm:disc} below.  

Suppose a reductive subgroup $L$ acts properly 
 and transitively on $X=G/H$.  
Then $L_H := L \cap\nolinebreak H$ is compact.  
We may assume
 that $L_K := L \cap K$ is a maximal compact subgroup of $L$
 containing $L_H$, 
 after possibly replacing $L$ by some conjugate.
Then the locally {\it{pseudo-Riemannian}} symmetric space
 $X=\Gamma \backslash G/H$ fibers 
 over the {\it{Riemannian}} locally symmetric space
 $Y_{\Gamma}= \Gamma \backslash L/L_K$ with fiber $F := L_K/L_H$:
\begin{equation}
\label{eqn:fiber}
  F \longrightarrow X_{\Gamma} \longrightarrow Y_{\Gamma}.  
\end{equation}

To expand functions on $X_{\Gamma}$ along the fiber $F$, 
 we define an endomorphism $p_{\tau}$
 of $C^{\infty}(X_{\Gamma})$ by
\[
  (p_{\tau} f)(\cdot)
  :=
  \frac{1}{\dim \tau}\int_K f(\cdot\,k)\operatorname{Trace} \tau(k)\ {\mathrm{d}} k 
\]
for every  $\tau \in \widehat{L_K}$. 
Then $p_{\tau}$ is an idempotent,
 namely,
 $p_{\tau}^2=p_{\tau}$.  
The $\tau$-component of $C^{\infty}(X_{\Gamma})$ 
 is defined by 
\[
   C^{\infty}(X_{\Gamma})_{\tau}
   :=
   \operatorname{Image}(p_{\tau})
   =\operatorname{Ker}(p_{\tau}-\operatorname{id}).  
\]
We note that
 $C^{\infty}(X_{\Gamma})_{\tau} \ne \{0\}$
 if and only if $\tau$
 has a nonzero $L_H$-invariant vector,
 i.e.\ $\tau \in (\widehat {L_K})_{L_H}$.  
It is easy to see
 that the projection $p_{\tau}$ commutes
 with any element in $\mathcal{Q}$ $(\simeq {\mathrm{d}}\ell({\mathfrak{Z}}({\mathfrak{l}}_{\mathbb{C}})))$, 
 but not always with \lq\lq{intrinsic differential operators}\rq\rq\ $D_{\Gamma} \in \mathcal{P}$
 $(\simeq {\mathbb{D}}_G(X))$, 
 and consequently it may well happen that
\[
   p_{\tau}(C^{\infty}(X_{\Gamma};{\mathcal{M}}_{\lambda}))
 \not \subset 
  C^{\infty}(X_{\Gamma};{\mathcal{M}}_{\lambda}).  
\]

To make a connection between the two subalgebras $\mathcal{P}$ and $\mathcal{Q}$, 
 we introduce a third subalgebra $\mathcal{R}$ of ${\mathbb{D}}(X_\Gamma)$,
 coming from the fiber $F$ in \eqref{eqn:fiber}.  
Namely, $\mathcal{R}$ is isomorphic to the ${\mathbb{C}}$-algebra ${\mathbb{D}}_{L_K}(F)$ of $L_K$-invariant differential operators $D$ on $F$, and obtained by extending elements of ${\mathbb{D}}_{L_K}(F)$
 to $L$-invariant differential operators on $X$, 
 yielding differential operators
 on the quotient $X_{\Gamma}$.  

Suppose now that we are in the setting \ref{setting:main}.  
The subgroup $L$ acts transitively on $X$
 by \cite[Lem.\,4.2]{KO13}
 and \cite[Lem.\,5.1]{kob94}.  
Moreover, 
 we can prove that
\begin{equation}
\label{eqn:QPR}
\mathcal{Q} \subset \langle \mathcal{P}, \mathcal{R}\rangle
\end{equation}
where $\langle \mathcal{P}, \mathcal{R} \rangle$ denotes 
 the subalgebra of ${\mathbb{D}}(X_{\Gamma})$
 generated by $\mathcal{P}$ and $\mathcal{R}$.  
This implies the following strong constraints
 on the restriction of representations:
 
\begin{theorem} \label{thm:disc}
In the setting \ref{setting:main}, 
 any irreducible $({\mathfrak{g}},K)$-module
 occurring in $C^{\infty}(X)$
 is discretely decomposable
 as an $({\mathfrak{l}},L \cap K)$-module.
\end{theorem}

See \cite{kob94, kob98b, kob98c}
 for a general theory
 of discretely decomposable restrictions of representations.
See also \cite{kob17}
 for a discussion on Theorem~\ref{thm:disc}
 when dropping the assumption
 that $L$ acts properly on $X$.  

In addition to \eqref{eqn:QPR}, 
 the quotient fields of $\mathcal{P}$ and $\langle \mathcal{Q}, \mathcal{R}\rangle$
 coincide \cite{kkdiffop}, 
 and we obtain the following:

\begin{theorem}[Transfer map]
\label{thm:transfer}
In the setting~\ref{setting:main}, for any $\tau \in (\widehat{L_K})_{L_H}$
 there is an injective map 
\[
  {\boldsymbol{\nu}}(\cdot, \tau)\colon
  \operatorname{Hom}_{\mathbb{C}\operatorname{-alg}}
  ({\mathbb{D}}_G(X), {\mathbb{C}})
  \hookrightarrow
  \operatorname{Hom}_{\mathbb{C}\operatorname{-alg}}
  ({\mathfrak{Z}}({\mathfrak{l}}_{\mathbb{C}}), {\mathbb{C}})
\]
such that for any $\lambda\in\operatorname{Hom}_{\mathbb{C}\operatorname{-alg}}
  ({\mathbb{D}}_G(X), {\mathbb{C}})$, any $f\in C^{\infty}(X_{\Gamma};{\mathcal{M}}_{\lambda})$, and any $z\in {\mathfrak{Z}}({\mathfrak{l}}_{\mathbb{C}})$,
\[
  {\mathrm{d}}\ell(z)(p_{\tau}f)
=
{\boldsymbol{\nu}}(\lambda,\tau)(z) \ p_{\tau}f.
\]
\end{theorem}

We write ${\boldsymbol{\lambda}}(\cdot, \tau)$ for the inverse map
 of ${\boldsymbol{\nu}}(\cdot,\tau)$ on its image.  
We call ${\boldsymbol{\nu}}
(\cdot,\tau)$ and ${\boldsymbol{\lambda}}(\cdot,\tau)$
 {\it{transfer maps}},
 as they \lq\lq{transfer}\rq\rq\
 eigenfunctions for $\mathcal{P}$ 
 to those for $\mathcal{Q}$,
 and vice versa, 
 on the $\tau$-component $C^{\infty}(X_{\Gamma})_{\tau}$.  

For an explicit description of transfer maps, 
 let 
\[
   \Phi^{\ast}\colon
   \operatorname{Hom}_{\mathbb{C}\text{-}\mathrm{alg}}
({\mathfrak{Z}}({\mathfrak{l}}_{\mathbb{C}}),\mathbb{C}) 
\overset \sim \to
 {\mathfrak{t}}_{\mathbb{C}}^{\ast}/W({\mathfrak{l}}_{\mathbb{C}})
\]
be the Harish-Chandra isomorphism
 as in \eqref{eqn:HCisom}, 
where $W({\mathfrak{l}}_{\mathbb{C}})$ denotes the Weyl group
 of the root system $\Delta({\mathfrak{l}}_{\mathbb{C}}, {\mathfrak{t}}_{\mathbb{C}})$
 with respect to a Cartan subalgebra ${\mathfrak{t}}_{\mathbb{C}}$
 in ${\mathfrak{l}}_{\mathbb{C}}$.  
We note that there is no natural inclusion relation
 between ${\mathfrak{j}}_{\mathbb{C}}$ and ${\mathfrak{t}}_{\mathbb{C}}$.  

For each $\tau \in (\widehat {L_K})_{L_H}$, 
we find an affine map
 $S_{\tau} \colon {\mathfrak{j}}_{\mathbb{C}}^{\ast} \to {\mathfrak{t}}_{\mathbb{C}}^{\ast}$
 such that the following diagram commutes:
$$
\xymatrix
{
{\mathfrak{j}}_{\mathbb{C}}^{\ast} 
\ar[d] 
\ar[r]^{S_{\tau}}
& {\mathfrak{t}}_{\mathbb{C}}^{\ast} \ar[d]
\\
{\mathfrak{j}}_{\mathbb{C}}^*/W 
& \mathfrak{t}_{\mathbb{C}}^{\ast}/W({\mathfrak{l}}_{\mathbb{C}})
\\
\operatorname{Hom}_{{\mathbb{C}}\text{-}\mathrm{alg}}({\mathbb{D}}_G(X),{\mathbb{C}}) \ar[u]^{\Psi^{\ast}}_{\text{\rotatebox{90}{$\sim$}}} \ar@{-->}[r]^{{\boldsymbol{\nu}}(\cdot,\tau)}
& 
\operatorname{Hom}
_{\mathbb{C}\text{-}\mathrm{alg}}
({\mathfrak{Z}}({\mathfrak{l}}_{\mathbb{C}}),\mathbb{C}) 
\ar[u]_{\Phi^{\ast}}^{\text{\rotatebox{-90}{$\sim$}}}
}
$$

Then a closed formula for the transfer map ${\boldsymbol{\nu}}(\cdot,\tau)$
 is derived from that of the affine map $S_{\tau}$
 which was determined explicitly 
in \cite{kkdiffop} for
the complexifications of the triples $(L,G,H)$
 in the setting \ref{setting:main}. 

Via the transfer maps,
 we can utilize representations
 of the subgroup $L$ 
 efficiently for the spectral analysis on $X_{\Gamma}$,
 as follows.  
As in \eqref{eqn:direct}, 
 let 
\begin{equation}
\label{eqn:GammaL}
   L^2(\Gamma \backslash L)
   \simeq
  \int_{\widehat L}^{\oplus}
   m_{\Gamma} (\vartheta) \, \vartheta \, \mathrm{d} \sigma(\vartheta)
\end{equation}
 be a disintegration of 
 the regular representation $L^2(\Gamma \backslash L)$
 of the subgroup $L$.  
Then the transform ${\mathcal{F}}_{\lambda}$ 
 in Theorem~\ref{thm:Fourier} can be built
 naturally by using \eqref{eqn:GammaL}
 and the expansion of $C_c^{\infty}(X_{\Gamma})$
 along the fiber $F$ in \eqref{eqn:fiber}.  
Consider the map
\begin{align*}
  \Lambda \colon
  (\widehat L)_{L_H} \times (\widehat{L_K})_{L_H} 
  &\to 
  \operatorname{Hom}_{{\mathbb{C}}\text{-}\mathrm{alg}}({\mathbb{D}}_G(X),{\mathbb{C}}), 
\end{align*}
$(\vartheta, \tau)\mapsto {\boldsymbol{\lambda}}(\chi_{\vartheta}, \tau)$, where $\chi_{\vartheta} \in 
\operatorname{Hom}
_{\mathbb{C}\text{-}\mathrm{alg}}
({\mathfrak{Z}}({\mathfrak{l}}_{\mathbb{C}}),\mathbb{C})$
 is the infinitesimal character of $\vartheta \in \widehat L$.  
Then the Plancherel measure $\mathrm{d} \mu$ on $\operatorname{Hom}_{{\mathbb{C}}\text{-}\mathrm{alg}}({\mathbb{D}}_G(X),{\mathbb{C}})$ in Theorem~\ref{thm:Fourier} can be
 defined by 
\begin{equation}
\label{eqn:dmu}
\mathrm{d} \mu= \Lambda_{\ast}(\mathrm{d} \sigma|_{(\widehat L)_{L_H}} \times (\widehat {L_K})_{L_H})
\end{equation}  

Detailed proofs of Theorems~\ref{thm:Fourier}, \ref{thm:selfadj}, 
 \ref{thm:mainII}, \ref{thm:disc}, and \ref{thm:transfer} will appear elsewhere.


\end{document}